\documentclass[twoside]{article}

\usepackage[accepted]{aistats2024}
\usepackage[round]{natbib}

\usepackage[utf8]{inputenc} 
\usepackage[T1]{fontenc}    
\usepackage{hyperref}       
\usepackage{url}            
\usepackage{booktabs}       
\usepackage{amsfonts}       
\usepackage{nicefrac}       
\usepackage{microtype}      
\usepackage{xcolor}         

\usepackage{amsmath}
\usepackage{amssymb}
\usepackage{mathtools}
\usepackage{amsthm}
\usepackage{enumitem}
\usepackage{booktabs}
\usepackage{mathtools}				
\DeclareMathOperator{\rk}{rank}	
\DeclareMathOperator{\tr}{trace}	    
\DeclareMathOperator{\vecc}{vec}	    

\newcommand{\RR}{\mathbb R}

\DeclarePairedDelimiter{\norm}{\lVert}{\rVert}

\usepackage{amsthm}
\newtheorem{theorem}{Theorem}

\newtheorem{lemma}{Lemma}
\theoremstyle{definition}
\newtheorem{definition}{Definition}

\theoremstyle{remark}

\usepackage{subcaption}

\begin{document}

%
\runningtitle{Absence of spurious solutions far from ground truth}

%
\runningauthor{Ma, Chen, Lavaei, Sojoudi}

\twocolumn[

\aistatstitle{Absence of spurious solutions far from ground truth: \\ A low-rank analysis with high-order losses}

\aistatsauthor{ Ziye Ma$^{*,1}$ \And Ying Chen$^{*,2}$ \And Javad Lavaei$^2$ \And Somayeh Sojoudi$^{1}$}

\aistatsaddress{Department of Electrical Engineering and Computer Science$^1$\\
   Department of Industrial Engineering and Operations Research$^2$ \\
   University of California, Berkeley } 
   
   ]

\begin{abstract}

Matrix sensing problems exhibit pervasive non-convexity, plaguing optimization with a proliferation of suboptimal spurious solutions. Avoiding convergence to these critical points poses a major challenge. This work provides new theoretical insights that help demystify the intricacies of the non-convex landscape.  In this work, we prove that under certain conditions, critical points sufficiently distant from the ground truth matrix exhibit favorable geometry by being strict saddle points rather than troublesome local minima. Moreover, we introduce the notion of higher-order losses for the matrix sensing problem and show that the incorporation of such losses into the objective function amplifies the negative curvature around those distant critical points. This implies that increasing the complexity of the objective function via high-order losses accelerates the escape from such critical points and acts as a desirable alternative to increasing the complexity of the optimization problem via over-parametrization. By elucidating key characteristics of the non-convex optimization landscape, this work makes progress towards a comprehensive framework for tackling broader machine learning objectives plagued by non-convexity.
\end{abstract}

\section{INTRODUCTION}
\label{sec:intro}
The optimization landscape of non-convex problems is notoriously complex to analyze in general due to the existence of an arbitrary number of spurious solutions (a spurious solution is a second-order critical point that is not a global minimum). As a result, if a numerical algorithm is not initialized close enough to a desirable solution, it may converge to one of those problematic spurious solutions. It may be acceptable (depending on the application) if the algorithm finds a critical point different from but close to the true solution, while converging to a point faraway implies the failure of the algorithm. In this paper, we study this issue by focusing on a class of benchmark non-convex problems, named matrix sensing, and analyze the landscape of the optimization problem in areas far away from the ground truth. 

To be more concrete, we focus on the Burer-Monteiro (BM) form of the problem:
\begin{equation}\label{eq:main}
\begin{aligned}
	\min_{X \in \mathbb{R}^{n \times r}} f(X) &\coloneqq \frac{1}{2} \|\mathcal{A}(XX^T) - b\|^2 
\end{aligned}
\end{equation}
where $b = \mathcal{A}(M^*) \in \RR^m$ is the vector of observed measurements and $M^*\in \RR^{ n \times n}$ is the ground truth solution to be found. We focus on the case where $M^*$ is positive semidefinite and symmetric since the asymmetric case (with $M^*$ being sign indefinite or rectangular) can be equivalently converted to the symmetric case \citep{bi2021local}. The linear operator $\mathcal{A}(\cdot)$ in \eqref{eq:main} is defined as $\mathcal{A}(M)= [ \langle A_1, M \rangle, \dots, \langle A_m, M \rangle]^T$, where $\{ A_i \}_{i=1}^m$ are $m$ sensing matrices, which can be assumed to be symmetric without loss of generality \citep{zhang2021general}. We use $r^*$ to denote the true rank of $M^*$. In \eqref{eq:main}, $X$ can have a search rank of $r$, which has to satisfy $r \geq r^*$. Since $r$ is often significantly smaller than $n$ in practice, the above factorization form optimizing over the matrix $X$ with $nr$ entries rather than a matrix $M$ with $n^2$ entries has major computational advantages.  

The matrix sensing problem is a canonical problem bearing many important applications, such as the matrix completion problem/netflix problem \citep{candes2009exact, candes2010power}, the compressed sensing problem \citep{donoho2006compressed}, the training of quadratic neural networks \citep{li2018algorithmic}, and an array of localization/estimation problems \citep{zhang2017conic, jin2019towards,singer2011angular, boumal2016nonconvex, shechtman2015phase,fattahi2020exact}. As a result, a better understanding of \eqref{eq:main} not only helps with the above applications, but also paves the way for the analysis of a broader range of non-convex problems. This is due in part to the fact that any polynomial optimization can be converted into a series of matrix sensing problems under benign assumptions \citep{molybog2020conic}. 

The major drawback of \eqref{eq:main} is that it may have an arbitrary number of spurious solutions, which cause ubiquitous local search algorithms to potentially end up with unwanted solution. Therefore, there has been an extensive investigation of the non-convex optimization landscape of \eqref{eq:main}, and the centerpiece notion is the restricted isometry property (RIP), defined below.
\begin{definition}[RIP] \citep{candes2009exact} \label{def:rip}
    Given a natural number $p$, the linear map $\mathcal{A}: \mathbb{R}^{n \times n} \mapsto \mathbb{R}^{m}$ is said to satisfy $\delta_{p}$-RIP if there is a constant $\delta_{p} \in [0,1)$ such that
    \[ (1-\delta_{p})\|M \|_F^2 \leq \| \mathcal{A} (M) \|^2 \leq (1 + \delta_{p}) \|M\|_F^2 \]
    holds for matrices $M \in \mathbb{R}^{n \times n}$ satisfying $\rk(M) \leq p$.
\end{definition}

Intuitively speaking, a smaller RIP constant means that the problem is easier to solve. For instance, if $\delta_{2r^*} = 0$, then $\mathcal{A}(\cdot)$ becomes the identity operator with $b = \vecc(M^*)$, which makes the problem trivial to solve for $M^*$. 

In this section, we provide a review of how RIP plays a central role in determining the optimization landscape of the non-convex problem \eqref{eq:main}, and explain why this problem still needs a further investigation even with the abundance of literature dedicated to this topic. To further streamline our presentation, we divide our discussion into two parts: RIP being smaller than $1/2$ and RIP being greater than $1/2$. 

\subsection{When RIP constant is smaller than $1/2$}

The attention to the RIP constant was first popularized by the study of using a convex semidefinite programming (SDP) relaxation to solve the matrix sensing problem \citep{recht2010guaranteed,candes2010power}. It was proven that as along as $\delta_{5r^*} \leq 1/10$, the SDP relaxation was tight and $M^*$ could be recovered exactly. Subsequently, \cite{bhojanapalli2016global} analyzed the factorized problem \eqref{eq:main} and concluded that as long as $\delta_{2r} \leq 1/5$, all second-order critical points (SOPs) of \eqref{eq:main} are ground truth solutions. \cite{zhu2018global,li2019non} also proved that $\delta_{4r} \leq 1/5$ is sufficient for the global recovery of $M^*$ under an arbitrary objective function (instead of the least-squares one in \eqref{eq:main}).  Later, by using a "certification of in-existence" technique, \cite{zhang2019sharp} established that $\delta_{2r} = 1/2$ was a sharp bound when $r=r^*$, meaning that as long as  $\delta_{2r} < 1/2$, all problem instances of \eqref{eq:main} are free of spurious solutions, and once $\delta_{2r} \geq 1/2$, it is possible to establish counter-examples with SOPs not corresponding to ground truth solutions. This aforementioned approach is important because it quantifies how restrictive RIP needs to be in order to ensure a benign landscape. 

Following the above line of work, \cite{ma2022sharp} proved that the same RIP bound of $1/2$ is also applicable in noisy cases, and \cite{bi2021local} showed that even for general low-rank optimization problems beyond matrix sensing, $\delta_{2r} < 1/2$ is still a sharp bound for the global recoverability of $M^*$, highlighting the importance of the $1/2$ bound.

Furthermore, when the RIP constant is small enough, various desirable properties hold, including fast convergence \citep{li2018algorithmic, wang2017unified} and spectral contraction \citep{stoger2021small,jin2023understanding}.

\subsection{When RIP constant is larger than $1/2$}

As proven in \cite{zhang2021general}, when the RIP constant of the problem is larger or equal to $1/2$, counterexamples can be found for which some SOPs are not global solutions. This means that in the regime where $\delta \geq 1/2$, the optimization landscape of \eqref{eq:main} becomes complex. Several works have attempted to provide limited mathematical guarantees in that regime.

\textbf{Benign landscape near $M^*$}. \cite{zhang2019sharp} proved that when $\delta_{2r} \geq 1/2$ for $r=1$, we can ensure the absence of spurious solutions in a local region that is close to $M^*$, depending on the RIP constant and also the size of $M^*$. \cite{zhang2020many} expanded that analysis to the regime of general $r$. Subsequently, \cite{ma2023noisy} extended the result to noisy and general objectives, proving the ubiquity of this phenomenon.

\textbf{Over-parametrization with $r \geq r^*$.} This line of work is concerned with the case where the search rank $r$ is greater than the true rank, which means that the complexity of the algorithm is increased. \cite{zhang2022improved} proved that if $r > r^*[(1+\delta_n)/(1-\delta_n)-1]^2/4$, with $r^* \leq r < n$, then every SOP $\hat X$ satisfies that $\hat X \hat X^\top = M^*$. \cite{ma2022global} derived a similar result for $l_1$ loss under an RIP-type condition.

\textbf{The SDP approach.} This approach uses the conventional technique of convex relaxations to solve the matrix sensing problem \citep{recht2010guaranteed,candes2011tight}. It was recently proven in \cite{yalcin2022semidefinite} that as long as the RIP constant $\delta_{2r^*}$ is lower than the maximum of $1/2$ and $2r^*/(n+(n-2r^*)(2l-5))$, the global solution of the SDP relaxation corresponds to $M^*$. Since this bound approaches 1 as $r^*$ increases, it is more appealing than using the factorized version \eqref{eq:main} when the RIP constant is not small.

\textbf{Lifting into tensor space.} Recently, \cite{ma2023over} proposed to lift the matrix decision variables of \eqref{eq:main} into tensors and then optimize over tensors. This approach is inspired by the Sum-of-Squares hierarchies and the authors proved that spurious solutions will be converted into saddle points through lifting, further alleviating the highly non-convex landscape of \eqref{eq:main} when $\delta_{2r}$ is close to 1. However, the shortcomings of this approach is that the complexity of the problem will increase exponentially if a high-order lifting is required. 

Overall, although various studies have been conducted to address the optimization landscape of \eqref{eq:main} when the RIP constant is larger than $1/2$, they either require to increase the complexity of the algorithm by a large margin (via over-parametrization  $r \gg r^*$, SDP relaxation, or tensor optimization) or require to initialize the algorithm close to  $M^*$. Therefore, the following question arises: \textit{Does there exist meaningful global guarantees for \eqref{eq:main} in the case of $\delta \geq 1/2$ without increasing the computational complexity of the problem drastically? }. In this paper, we offer a partial affirmative answer to this question via our notion of high-order losses. 

\subsection{Main Contributions}
\begin{enumerate}
	\item We prove that all critical points of \eqref{eq:main} are strict saddles when reasonably away from $M^*$, with the intensity of the smallest eigenvalue of the Hessian at each saddle point being proportional to its distance to $M^*$. This result implies that there are no spurious solutions far away from $M^*$, and that it is possible to reach a vicinity of $M^*$ with saddle-escaping algorithms even with poor initializations.
	\item As a by-product of the above result, we derive sufficient conditions on $M^*$ to ensure that there are no spurious solutions in the entire space even in the regime of high RIP constants. 
	\item We introduce the notion of high-order losses where a penalization term with a controllable degree is added to the objective function of \eqref{eq:main} with the property that the penalty is zero at the ground truth. We show that critical points far away from $M^*$ will still remain strict saddles, while the spurious solutions will be easier to escape as the degree of the loss increases. In other words, the landscape of the optimization problem is reshaped favorably by the inclusion of such penalties. Our result implies that increasing the complexity of the objective function serves as a viable alternative to increasing the complexity of the problem via over-parametrization. 
\end{enumerate}

\subsection{Notations}
The notation $M \succeq 0$ means that $M$ is a symmetric and positive semidefinite (PSD) matrix.  $\sigma_i(M)$ denotes the $i$-th largest singular value of a matrix $M$, and $\lambda_i(M)$ denotes the $i$-th largest eigenvalue of $M$. $\lambda_{\min}(M)$ and $\lambda_{\max}(M)$ respectively denotes the minimum and maximum eigenvalues of $M$. $\norm{v}$ denotes the Euclidean norm of a vector $v$, while $\norm{M}_F$ and $\norm{M}_p$ denote the Frobenius norm and induced $l_p$ norm of a matrix $M$, respectively for $p \geq 2$. $\langle A,B \rangle$ is defined to be $\tr(A^TB)$ for two matrices $A$ and $B$ of the same size. For a matrix $M$, $\vecc(M)$ is the usual vectorization operation by stacking the columns of the matrix $M$ into a vector.  $[n]$ denotes the integer set $\{1,\dots,n\}$.

The Hessian of the function $f(X)$ in \eqref{eq:main}, denoted as $\nabla^2 f(X): \RR^{n \times r} \times \RR^{n \times r} \mapsto \RR$, can be regarded as a quadratic form whose action on any two matrices $U,V \in \RR^{n \times r}$ is given by
\[
\nabla^2 f(X)[U,V]=\sum_{i,j,k,l=1}^{n,r,n,r}\frac{\partial^2}{\partial X_{ij}\partial X_{kl}}f(X)U_{ij}V_{kl}.
\]

\section{DISAPPEARANCE OF SPURIOUS SOLUTIONS FAR FROM GROUND TRUTH}
\label{sec:l2_case}

As discussed in Section \ref{sec:intro}, the optimization landscape of \eqref{eq:main} is benign (in the sense of having no spurious solutions) if $\delta_{2r} < 1/2$ and benign in a region close to $M^*$ if $\delta_{2r} \geq 1/2$. In this section, we study the landscape far away from $M^*$ in the problematic case $\delta_{2r} \geq 1/2$. To do so, we focus on the first-order critical points, and study the eigenvalues of the Hessian at these points because if they exhibit negative eigenvalues, it means that these first-order critical points are strict saddles, possessing escape directions. Before diving into more details, we present the first- and second-order optimality conditions for \eqref{eq:main}:

\begin{lemma}\label{lem:l2_crit_cond}
	A point $X$ is a first-order critical point of \eqref{eq:main} if 
	\begin{equation}
		\nabla f(X) = \left(\sum_{i=1}^m \langle A_i, XX^\top - M^* \rangle  A_i \right) X = 0 \label{eq:focp_l2}
	\end{equation}
	and it is a second-order critical point if it satisfies the above condition together with 
	\begin{equation}\label{eq:socp_l2}
		\begin{aligned}
			\nabla^2 f(X)[U,U] =  \sum_{i=1}^m 	 \langle A_i, UX^\top + XU^\top \rangle^2 + \\
		\langle A_i, XX^\top - M^* \rangle \langle A_i, 2 UU^\top \rangle  \geq 0 \quad \forall U \in \RR^{n \times r} 
		\end{aligned}
	\end{equation}
\end{lemma}

The proof of this lemma is plain calculus thus omitted for simplicity. Focusing on \eqref{eq:socp_l2}, it is apparent that for a first-order critical point $\hat X$ satisfying $\nabla f(\hat X) = 0$, the Hessian $\nabla^2 f(\hat X)[U,U]$ can be broken down into the summation of two terms:
\begin{align*}
	T_1 &\coloneqq \sum_{i=1}^m 	 \langle A_i, U\hat X^\top + \hat XU^\top \rangle^2 = \|\mathcal{A}(U\hat X^\top + \hat XU^\top)\|_2^2, \\
	T_2 &\coloneqq \sum_{i=1}^m  \langle A_i, \hat X \hat X^\top - M^* \rangle \langle A_i, 2 UU^\top \rangle \\
	&= 2 \langle \nabla h(\hat X \hat X^\top) , UU^\top \rangle
\end{align*}
Assuming that the problem \eqref{eq:main} satisfies the RIP condition with some constant $\delta_p$ for $p \geq 2r^*$, one can write  
\[
	T_1 \leq (1+\delta_p) \|U\hat X^\top + \hat XU^\top\|^2_F,
\]
which means that $T_1$ can be upper-bounded naturally since $U$ can be assumed to have a unit scale without loss of generality. Therefore, if we can somehow show that there exists $U \in \RR^{n \times r}$ to make $T_2$ negative with a sufficiently large magnitude, then $\hat X$ becomes a saddle point. Combining the RIP condition and mean value theorem, we know that 
\begin{align*}
	h(M^*) \geq &h(\hat X \hat X^\top) + \langle \nabla h(\hat X \hat X^\top), M^* - \hat X \hat X^\top \rangle \\
	+& \frac{1-\delta_p}{2} \|\hat X \hat X^\top -M^* \|^2_F
\end{align*}
where $h(\cdot)$ is defined as
\begin{equation}
	h(M) = \frac{1}{2} \|\mathcal{A}(M - M^*)\|^2
\end{equation}
Given $\nabla f(\hat X) = 0$ and the expression in \eqref{eq:focp_l2}, we obtain that
 \[
 	\langle \nabla h(\hat X \hat X^\top), M^* \rangle \leq - \frac{1-\delta_p}{2} \|\hat X \hat X^\top -M^* \|^2_F
 \]
 since $h(\hat X \hat X^\top) \geq h(M^*)$ by definition. This implies that there exist directions that make $T_2$ have large negative values when $\hat X \hat X^\top$ is far away from $M^*$. Expanding on this simple observation, we formally establish Theorem~\ref{thm:l2_nospurious}, serving as the cornerstone of all results in this paper. A detailed proof can be found in the Appendix.
 
\begin{theorem} \label{thm:l2_nospurious}
	Assume that \eqref{eq:main} satisfies the RIP$_{r+r^*}$ property with constant $\delta \in [0,1)$. Given a first-order critical point $\hat X \in 
	\RR^{n \times r}$ of \eqref{eq:main}, if it satisfies the inequality
	\begin{equation}\label{eq:l2_criterion}
		\|\hat X \hat X^\top -M^* \|^2_F > 2 \frac{1+\delta}{1-\delta} \tr(M^*) \sigma_r(\hat X)^2,
	\end{equation} 
	then $\hat X$ is not a second-order critical point and is a strict saddle point with $\nabla^2 f(\hat X)$ having a strictly negative eigenvalue not larger than 
	\begin{equation}\label{eq:hessian_eigenval_l2}
		2 (1+\delta) \sigma_r(\hat X)^2 - \frac{\|\hat X \hat X^\top -M^* \|^2_F (1-\delta)}{\tr(M^*)}
	\end{equation}
\end{theorem}
Theorem~\ref{thm:l2_nospurious} states that if a first-order critical point is far from the ground truth, it cannot be a spurious solution, and always exhibits an escape direction with its magnitude proportional to the squared distance between $\hat X \hat X^\top$ and $M^*$. This further elucidates the fact that even if \eqref{eq:main} is poorly initialized, it is possible to converge to a vicinity of $M^*$ with saddle-escaping algorithms, which we will numerically illustrate in Section~\ref{sec:experiments}.

In contrary to Theorem~\ref{thm:l2_nospurious}, the existing results in the literature state that there are no spurious solutions in a small neighborhood of $M^*$ \citep{zhang2020many,bi2020global,ma2023geometric}. We recall a well-known result in the literature regarding this property below.
\begin{theorem}[\cite{zhang2020many}] \label{thm:near_m_star}
	Assume that \eqref{eq:main} satisfies the RIP property with constant $\delta \in [0,1)$. Given an arbitrary constant $\tau \in (0,1-\delta^2)$, if a second-order critical point $\hat X \in \RR^{n \times r}$ of \eqref{eq:main} satisfies
	\begin{equation}\label{eq:close_criterion}
		\|\hat X \hat X^\top - M^*\|_F \leq \tau \lambda_{r^*}(M^*)
	\end{equation}
	then $\hat X$ corresponds to the ground truth solution.
\end{theorem}

Theorem~\ref{thm:near_m_star} serves as a classic result in which a good enough initialization (the best possible scenario is initializing at the ground truth) can lead to the global solution easily, which could be observed in a wide range of problems. However, to the contrary, Theorem~\ref{thm:l2_nospurious} shows that even with a very bad initialization, it is possible (via saddle-escaping methods) to reach a vicinity of the ground truth due to the lack of spurious solutions far away, which complements the prior wisdom nicely.

The existence of both Theorem~\ref{thm:l2_nospurious} and Theorem~\ref{thm:near_m_star} inspires us to think that it is possible for some matrix sensing problems to have no spurious solutions at all (since both far away and near regions of $M^*$ are benign). Therefore the question becomes whether there is a sufficient condition to ensure that the two regions \eqref{eq:l2_criterion} and \eqref{eq:close_criterion} overlap? An affirmative answer guarantees that the entire optimization landscape has no spurious solutions, even if $\delta > 1/2$, exceeding the sharp RIP bound regardless of $\mathcal{A}$ and $M^*$ \cite{zhang2019sharp}. In what follows, we will prove that if $\|M^*\|_F$ is small, the regions \eqref{eq:l2_criterion} and \eqref{eq:close_criterion} overlap.
\begin{theorem}\label{thm:small_mstar}
	Consider the problem \eqref{eq:main} under the RIP$_{r+r^*}$ property with a constant $\delta \in [0,1)$. Assume that its ground truth solution $M^*$ satisfies the following inequality
	\begin{equation}\label{eq:mstar_criterion}
		\|M^*\|_F \frac{\tr(M^*)}{\lambda_{r^*}(M^*)} \leq \frac{\sqrt{r}}{2\sqrt{2}} \sqrt{\frac{(1-\delta)^5}{(1+\delta)}},
	\end{equation}
	Then, every second-order critical point $\hat X$ of \eqref{eq:main} satisfies
	\[
		\hat X \hat X^\top = M^*
	\]
\end{theorem}
The proof can be found in the Appendix.

\begin{table*}[h]
\centering
\begin{tabular}{l l l l l l l}
\toprule
n & $\lambda$ & $\lambda_{\min}(\nabla^{2} f^{l}(\hat{X}))$ & $\lambda_{\max}(\nabla^{2} f^{l}(\hat{X}))$ &$\lambda_{\min }(\nabla^{2} f^{l}(X^{*}))$ & $\lambda_{\max}(\nabla^{2} f^{l}(X^*))$ \\
\midrule
3 & 0 & 1.821 & 3.642 & 2.18 &  4.36 \\
3 & 0.5 & 1.779 & 3.855 & 2.18 &  4.36\\
3 & 5 & 1.594 & 7.422 & 2.18 &  4.36\\
3 & 50 & 1.470 & 55.028 & 2.18 &  4.36\\
5 & 0 & 0.429 & 3.898 & 0.54 &  4.72\\
5 & 0.5 & 0.421 & 4.106 & 0.54 &  4.72\\
5 & 5 & 0.385 & 9.117 & 0.54 &  4.72\\
5 & 50 & 0.354 & 69.816 & 0.54 &  4.72 \\
7 & 0 & 0.516 & 3.642 & 0.72 &  5.08 \\
7 & 0.5 & 0.502 & 4.122 & 0.72 & 5.08\\
7 & 5 & 0.456 & 10.006 & 0.72 &  5.08\\
7 & 50 & 0.433 & 75.786 & 0.72 &  5.08 \\
9 & 0 & 0.609 & 3.930 & 0.90 & 5.44 \\
9 & 0.5 & 0.601 & 4.315 & 0.90 & 5.44\\
9 & 5 & 0.557 & 10.915 & 0.90 & 5.44 \\
9 & 50 & 0.528 & 84.002 & 0.90 & 5.44 \\
\bottomrule
\end{tabular}
\caption{The smallest eigenvalue of the Hessian at a spurious local minimum $\hat{X}$ and ground truth $X^*$, with $\epsilon = 0.3$ and additional high-order loss function $l=4$ (note that $\hat X$ is not too far from $X^*$ since Theorem \ref{thm:l2_nospurious} shows that there are no such spurious solutions). The problem satisfies the $RIP_{2r}$-property with $\delta = \frac{1-\varepsilon}{1+\varepsilon}=0.538>1/2$, and hence has spurious local minima.}
\label{tb:least-eigenvalue}
\end{table*}

\section{HIGHER-ORDER LOSS FUNCTIONS }\label{sec:modified_loss}
\begin{figure*}[h]
\centering
    \begin{minipage}{0.45\textwidth}
    \centering
    \begin{subfigure}{\textwidth}
      \includegraphics[width=\linewidth]{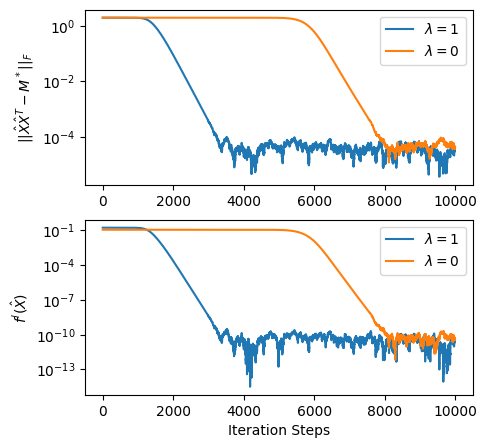}
      \caption{$\lambda = 0$ converges to ground truth}
    \end{subfigure}
  \end{minipage}
  \begin{minipage}{0.47\textwidth}
    \centering
    \begin{subfigure}{\textwidth}
      \includegraphics[width=\linewidth]{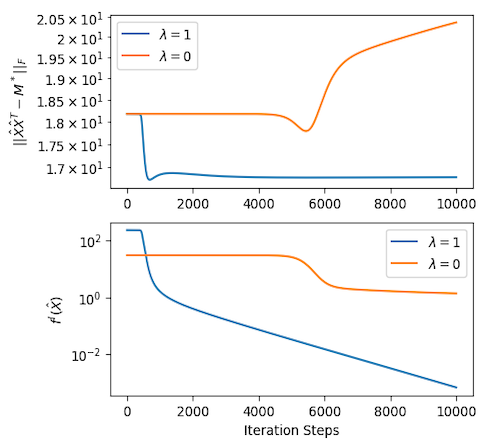}
      \caption{$\lambda = 0$ converges to a spurious solution around the ground truth}
    \end{subfigure}
  \end{minipage}%

    \caption{The evolution of the objective function and the error between the obtained solution $\hat{X}\hat{X}^T$ and the ground truth $M^*$ during the iterations of the perturbed gradient descent method, with a constant step-size. In both cases, high-order loss functions accelerate the convergence.}
  \label{fig:gaussian-pgd}
\end{figure*}
Although Theorem~\ref{thm:l2_nospurious} proves that critical points far away from the ground truth are strict saddle points, the time needed to escape such points depends on the local curvature of the function \citep{ge2017no,jin2021nonconvex}. Therefore, it is essential to understand whether the curvatures at saddle points could be enhanced to reshape the landscape favorably. In this section, we provide an affirmative answer to this question by using a modified loss function.

Our main goal in matrix sensing is to recover the ground truth matrix $M^*$ via $m$ measurements, and we minimize a mismatch error in \eqref{eq:main} to achieve this goal. An $l_2$ loss function is used in \eqref{eq:main} due to its smooth and nonnegative properties, which is the most common objective function in the machine learning literature. However, in this work, we introduce a high-order loss function as penalization, namely an $l_p$ loss function with $p > 2$, and show that this will reshape the landscape of the optimization problem. To be concrete, we propose to optimize over this modified problem:
\begin{equation}\label{eq:l_p_problem}
\begin{aligned}
	\min_{X \in \mathbb{R}^{n \times r}} f_{\lambda}^l(X) &\coloneqq  f(X) + \lambda f^l(X)
\end{aligned}
\end{equation}
where
\begin{subequations}
\begin{align}
f^l(X) &\coloneqq  \frac{1}{l} \|\mathcal{A}(XX^\top)-b\|^l_l\label{eq:lp_f_func}
\\
h^l(M) &\coloneqq  \frac{1}{l} \|\mathcal{A}(M) -b\|^l_l 
\label{eq:lp_h_func}
\end{align}
\end{subequations}

where $l \geq 2$ is an even natural number to ensure the non-negativity of the loss function and $\lambda > 0$ is a penalty coefficient. The intuition behind using a high-order objective can be easily demonstrated via the scalar example: 
\begin{equation}
	\min_{x \in \mathbb{R}} g(x) \coloneqq \frac{1}{l}(x^2-a)^l
\end{equation} 
for some constant $a \in \RR$ and an even number $l \geq 2$. This problem is a scalar analogy of $f^l(X)$ with $\mathcal{A}(\cdot)$ being the identity operator. In this example, the derivatives are 
\begin{align*}
	g'(x) &= 2x(x^2-a)^{l-1}, \\
	g''(x) &= 2 (x^2-a)^{l-2} \left[ (l-1) 2x^2  + (x^2-a) \right]
\end{align*}
It can be observed that as $l$ increases, the first- and second-order derivatives will be amplified, provided that $(x^2-a)$ is larger than one (i.e., our point is reasonably distant from the ground truth $a$). However, there is an issue with optimizing $g(x)$ directly, and we need to use $f_{\lambda}^l(X)$ instead of $f^l(X)$. If we directly minimize  $f^l(X)$ with $l > 2$, the Hessian at any point $X$ with the property $XX^T=M^*$ becomes zero, which makes the convergence  extremely slow as approaching the ground truth (with a sub-linear rate).  This is because the local convergence rate of descent numerical algorithms depends on the condition number of the Hessians around the solution \citep{wright2022optimization}. Conversely, when using the original objective \eqref{eq:main}, we see from \eqref{eq:socp_l2} that even if $XX^\top$ is close to $M^*$, the Hessian is positive semidefinite, and therefore adding $f^l(X)$ to the objective of \eqref{eq:main} will not change the sign of the Hessian around the solution. 

Secondly, if $l$ is large and $\|XX^\top-M^*\|_F$ is less than one, the term $\langle A_i, XX^\top - M^* \rangle^{l-1}$ appearing in the gradient of $f^l(X)$ (see Lemma~\ref{lem:lp_crit_cond} in Appendix) is very small due to its exponentiation nature. This means that minimizing $f^l(X)$ alone will suffer from the vanishing issue and slow growth rate in a local region around $M^*$. 

Due to the above reasons, we mix $f^l(X)$ with the original objective in \eqref{eq:main} and use the parameter $\lambda$ to control the effect of the penalty term, in an effort to balance local rate of convergence to $M^*$ and prominent eigenvalues of the Hessian at points far away from $M^*$. By using \eqref{eq:l_p_problem}, we can arrive at a similar result to Theorem~\ref{thm:l2_nospurious}
\begin{theorem} \label{thm:lp_nospurious}
	Assume that the operator $\mathcal{A}(\cdot)$ satisfies the RIP$_{r+r^*}$ property with constant $\delta \in [0,1)$. Consider the high-order optimization problem  \eqref{eq:l_p_problem} such that $l\geq2$ is even. Given a first-order critical point $\hat X \in 
		\RR^{n \times r}$ of \eqref{eq:l_p_problem}, if 
	\begin{equation}\label{eq:lp_criterion}
		D^2 \geq \tr(M^*) \sigma_r^2(\hat X) \frac{(1+\delta) + \lambda (l-1) (1+\delta)^{l/2} D^{l-2}}{(1-\delta)/2 + \lambda C(l) (1-\delta)^{l/2} D^{l-2}},
	\end{equation} 
	then $\hat X$ is a strict saddle point with $\nabla^2 f(\hat X)$ having a strictly negative eigenvalue not larger than  
	\begin{equation}\label{eq:hessian_eigenval_lp}
	\begin{aligned}
		& \left[ 2 (1+\delta) \sigma_r(\hat X)^2 - \frac{D^2 (1-\delta)}{\tr(M^*)} \right] + \\
		\lambda D^{l-2} & \left[ 2 (1+\delta)^{l/2} (l-1) \sigma_r(\hat X)^2 - 2 \frac{(1-\delta)^{l/2} C(l) D^2}{\tr(M^*)}\right]
	\end{aligned}
	\end{equation}
	where
	\begin{equation}\label{eq:D_C_def}
	\begin{aligned}
		D &\coloneqq \|\hat X \hat X^\top - M^* \|_F, \\
		C(l) &\coloneqq m^{(2-l)/2}\left( \frac{2^l-1}{l}-1\right)
	\end{aligned}
	\end{equation}
\end{theorem}
Theorem~\ref{thm:lp_nospurious} serves as a direct generalization of Theorem~\ref{thm:l2_nospurious}, as it recovers the statements of Theorem~\ref{thm:l2_nospurious} when $l$ is set to 2 or $\lambda$ is set to 0. By comparing \eqref{eq:hessian_eigenval_lp} to \eqref{eq:hessian_eigenval_l2}, the bound on the smallest eigenvalue of the Hessian has an additional term that is amplified by $\|\hat X \hat X^\top - M^* \|_F^{l-2}$. As a result, $\hat X$ has a more pronounced escape direction in \eqref{eq:l_p_problem} compared to \eqref{eq:main} when $\|\hat X \hat X^\top - M^* \|_F$ is large. Concerning the tightness of the bounds in Theorem~\ref{thm:l2_nospurious} and Theorem~\ref{thm:lp_nospurious}, they depend on three factors: the RIP constant, the Frobenius norm of $M^*$, and the smallest non-zero singular value of $\hat X$. While the RIP constant indicates problem difficulty and is immutable, knowing that $\|M^*\|_F$ is small (as in the extreme case shown in \eqref{eq:mstar_criterion}), or that $\hat X$ has a minimal singular value (which is computable), suggests that the bounds can be very tight. This implies that if $\hat X$ is a second-order point, it will be very close to $M^*$. Recent studies \citep{stoger2021small,jin2023understanding} have shown that a small random initialization can lead to small $r^{th}$ eigenvalues during the optimization process, resulting in tight bounds for our new Theorems.

Theorem~\ref{thm:lp_nospurious} contrasts well with an existing approach that utilizes a lifting technique to eliminate spurious solutions. \cite{ma2023over} states that by lifting the search space to the regime of tensors, a higher degree of parametrization can amplify the negative curvature of Hessian. In contrast, Theorem~\ref{thm:lp_nospurious} offers similar benefits by using a more complex objective function. This means that without resorting to massive over-parametrization, similar results can be achieved via using a more complex loss function. Having said that, the technique presented in \citep{ma2023over} can amplify the negative curvature of those points $X$ that satisfy
\[
	\|XX^\top -M^*\|_F^2 \geq \frac{1+\delta}{1-\delta} \tr(M^*) \sigma_r^2(\hat X)
\]
where in comparison to \eqref{eq:lp_criterion} the multiplicative factor to $\tr(M^*) \sigma_r^2(\hat X)$ becomes
\begin{equation*}
	\frac{(1+\delta) + \lambda (l-1) (1+\delta)^{l/2} D^{l-2}}{(1-\delta)/2 + \lambda C(l) (1-\delta)^{l/2} D^{l-2}},
\end{equation*}
which is on the order of magnitude of
\begin{equation*}
	\mathcal{O}\left( l \left(\frac{\sqrt{m}}{2}\right)^l\left(\frac{1+\delta}{1-\delta}\right)^{l/2} \right),
\end{equation*}
making the region for which this amplification can be observed smaller if $l$ is large. This means that by utilizing a high-order loss, we can recover some of the desirable properties of an over-parametrized technique, but a gap still exists due to the smaller parametrization. Combining a high-order loss function and a modest level of parametrization is left as future work. 

%

\section{SIMULATION EXPERIMENTS}
\label{sec:experiments}

This section serves to provide numerical validation for the theoretical findings presented in this paper. We will begin by investigating the behavior of the Hessian matrix when utilizing high-order loss function. Subsequently, we will showcase the remarkable acceleration in escaping saddle points achieved by employing Perturbed Gradient Descent in conjunction with high-order loss functions compared to the standard optimization problem \eqref{eq:main}. Lastly, we will provide a comparative illustration of the landscape both with and without the incorporation of high-order loss functions.\footnote{The code to this section can be found at https://github.com/anonpapersbm/high\_order\_obj}.

\begin{figure}[h]
    \centering
        \begin{minipage}{0.49\textwidth}
    \includegraphics[width = \linewidth]{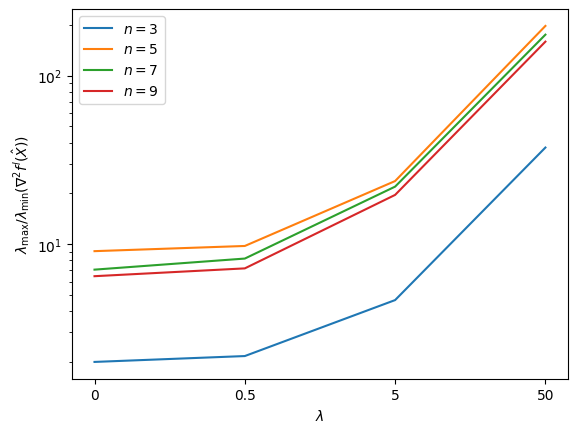}
      \end{minipage}
    \caption{The ratio between the largest and smallest eigenvalue of Hessian at the spurious local minimum $\lambda_{\max}/\lambda_{\min}(\nabla^{2} f^{l}(\hat{X}))$ with respect to $\lambda$ under different size $n$.}
        \label{fig:ratio}
\end{figure}

We first focus on a benchmark matrix sensing problem with the operator $\mathcal{A}$ defined as
\begin{equation}
\mathcal{A}_\epsilon(\mathbf{M})_{i j}:=\left\{\begin{array}{ll}
\mathbf{M}_{i j}, & \text { if }(i, j) \in \Omega \\
\epsilon \mathbf{M}_{i j}, & \text { otherwise }
\end{array},\right.
\label{eq:Yalcin-example}
\end{equation}

where $\Omega=\{(i, i),(i, 2 k),(2 k, i) \mid \forall i \in[n], k \in[\lfloor n / 2\rfloor]\}$, $0<\epsilon<1$. \citet{yalcin2022semidefinite} has proved that while satisfying RIP property with $\delta_{2 r}=(1-\epsilon) /(1+\epsilon)$, this problem has $\mathcal{O}\left(2^{\lceil n / 2\rceil}-2\right)$ spurious local minima. In order to analyze the influence of high-order loss functions on the optimization landscape, we conduct an analysis of the spurious local minima and of the ground truth matrix for both the vanilla problem \eqref{eq:main} and the altered problem \eqref{eq:l_p_problem} with $l=4$. We consider a spurious local minimum $\hat X$ (note that such points cannot be too far away from $M^*$ due to Theorem \ref{thm:l2_nospurious}).  The findings of this study are presented in Table \ref{tb:least-eigenvalue}, while the ratio between the largest and smallest eigenvalues at the spurious local minimum $\hat{X}$ is plotted in Figure~\ref{fig:ratio}.

\begin{figure*}[h]
\centering
  \begin{minipage}{\textwidth}
    \centering
    \begin{subfigure}{\textwidth}
      \includegraphics[width=\linewidth]{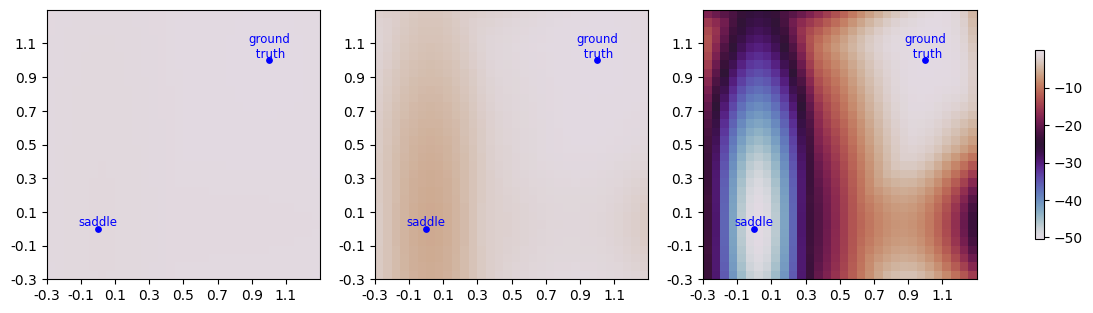}
    \end{subfigure}
  \end{minipage}%
  
\vspace{\baselineskip}
  \begin{minipage}{\textwidth}
    \centering
    \begin{subfigure}{\textwidth}
      \includegraphics[width=\linewidth]{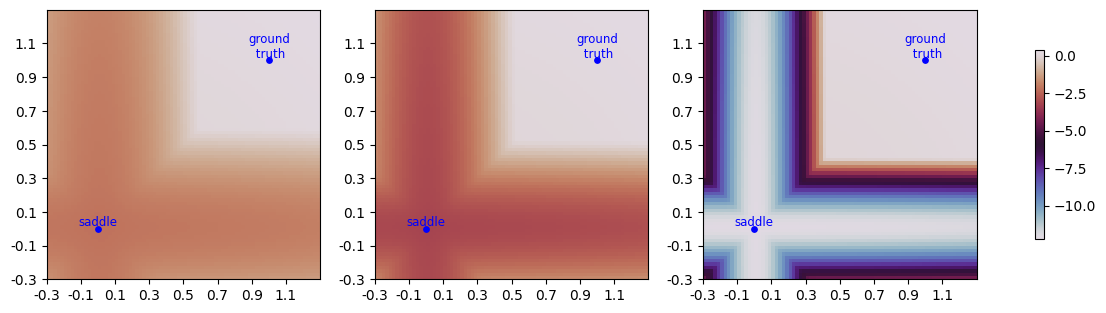}
  
    \end{subfigure}
  \end{minipage}%
    \caption{The value of the minimum eigenvalue of the Hessian around saddle points: The first row is for randomly generated Gaussian matrix with $m=20, n=20$, and the second row is for problem \eqref{eq:Yalcin-example} with $n =21, \epsilon = 0.1$. $\lambda = 0$ (left column), $\lambda = 0.5$ (middle column), $\lambda = 5$ (right column), with x-axis and y-axis as two orthogonal directions from the critical point to the ground truth.}
  \label{fig:landscape}
\end{figure*}

Table \ref{tb:least-eigenvalue} shows that as the intensity of high-order loss function increases via $\lambda$, the behavior of the Hessian eigenvalues exhibits distinct characteristics across different points in the optimization landscape. Specifically, the smallest and the largest eigenvalues of the Hessian at the ground truth matrix remain constant. In contrast, the smallest eigenvalue of the Hessian at the spurious local minimum, which is initially positive, decreases as $\lambda$ increases. This decreasing trend facilitates the differentiation of spurious local minima from the global minimum, as they become less favorable. Simultaneously, the largest eigenvalue of the Hessian matrix at the spurious local minimum increases at a significantly faster rate. This suggests that the incorporation of high-order loss functions amplifies the magnitude of the eigenvalues in the Hessian matrix at spurious local minima, increasing the ratio between the largest and smallest eigenvalues, while having no impact on those at the ground truth.

Following that, we will present the acceleration in effectively navigating away from saddle points. For randomly generated zero-mean Gaussian sensing matrices with i.i.d. entries, we apply small initialization and perturbed gradient descent which adds small Gaussian noise when the gradient is close to zero. In Figure \ref{fig:gaussian-pgd}, we compare the evolution of the distance from the ground truth matrix $\|\hat{X}\hat{X}^T-M^*\|_F$ and the value of the objective function $f^l(\hat{X})$. Although Figure \ref{fig:gaussian-pgd} demonstrates the behavior for a single problem, we observed the same phenomenon for many different trials. By incorporating a high-order loss function (specifically, with $l=4$), the optimization process exhibits enhanced convergence compared to the standard vanilla problem. This accelerated convergence can be attributed to the presence of a substantial negative eigenvalue of the Hessian matrix, which effectively facilitates the algorithm's escape from regions proximate to spurious local minimum.

Finally, we explore the optimization landscape in terms of the distance from the ground truth matrix and the intensity of high-order loss functions. We explore both random Gaussian sensing matrices with size $m=20, n =20$, and the problem \eqref{eq:Yalcin-example} with the parameters $n = 21$ and $l= 4$. The result is plotted in Figure \ref{fig:landscape}, where the x-axis and y-axis are two orthogonal directions from the critical point to the ground truth. By looking horizontally across the figure, we can observe that increasing the parameter $\lambda$ leads to the amplification of the least negative eigenvalue of the Hessian matrix at saddle points. As $\lambda$ increases, the least eigenvalue of the Hessian at this saddle point, which is initially negative, decreases further. This reduction in the magnitude of the negative eigenvalue makes it easier to escape from this saddle point during optimization. 

This example could also directly corroborate Theorem \ref{thm:lp_nospurious} (Thereby Theorem~\ref{thm:l2_nospurious}). For instance, when $\lambda = 0.5$, the minimum eigenvalue of Hessian matrix at the first-order critical point $\hat{X}$ is $\lambda_{\min}(\nabla^{2} f^{l}(\hat{X}))=-3.201$, which is smaller than the eigenvalue-bound $-2.274$; the distance from the ground truth matrix is $D \coloneqq \|\hat X \hat X^\top - M^* \|_F = 11.0$, larger than the distance bound in \eqref{eq:lp_criterion}, validating Theorem \ref{thm:lp_nospurious}.

\section{CONCLUSION}
This work theoretically establishes favorable geometric properties in those parts of the space far from the globally optimal solution for the non-convex matrix sensing problem. We introduce the notion of high-order loss functions and show that such losses reshape the optimization landscape and accelerate escaping saddle points. Our experiments demonstrate that high-order penalties decrease minimum Hessian eigenvalues at spurious points while intensifying ratios. Secondly, perturbed gradient descent exhibits accelerated saddle escape with the incorporation of high-order losses. Collectively, our theoretical and empirical results show that using a modified loss function could make non-convex functions easier to deal with and achieve some of the desirable properties of a lifted formulation without enlarging the search space of the problem exponentially.

\section{ACKNOWLEDGEMENTS}
This work was supported by grants from ARO, ONR, AFOSR, NSF, and the UC Noyce Initiative.
\bibliography{references.bib}

\begin{thebibliography}{}

\bibitem[Bhojanapalli et~al., 2016]{bhojanapalli2016global}
Bhojanapalli, S., Neyshabur, B., and Srebro, N. (2016).
\newblock Global optimality of local search for low rank matrix recovery.
\newblock In {\em Advances in Neural Information Processing Systems},
  volume~29.

\bibitem[Bi and Lavaei, 2020]{bi2020global}
Bi, Y. and Lavaei, J. (2020).
\newblock Global and local analyses of nonlinear low-rank matrix recovery
  problems.
\newblock arXiv:2010.04349.

\bibitem[Bi et~al., 2022]{bi2021local}
Bi, Y., Zhang, H., and Lavaei, J. (2022).
\newblock Local and global linear convergence of general low-rank matrix
  recovery problems.
\newblock {\em AAAI-22}.

\bibitem[Boumal, 2016]{boumal2016nonconvex}
Boumal, N. (2016).
\newblock Nonconvex phase synchronization.
\newblock {\em SIAM Journal on Optimization}, 26(4):2355--2377.

\bibitem[Candes and Plan, 2011]{candes2011tight}
Candes, E.~J. and Plan, Y. (2011).
\newblock Tight oracle inequalities for low-rank matrix recovery from a minimal
  number of noisy random measurements.

\bibitem[Cand{\`e}s and Recht, 2009]{candes2009exact}
Cand{\`e}s, E.~J. and Recht, B. (2009).
\newblock Exact matrix completion via convex optimization.
\newblock {\em Foundations of Computational Mathematics}, 9(6):717--772.

\bibitem[Cand{\`e}s and Tao, 2010]{candes2010power}
Cand{\`e}s, E.~J. and Tao, T. (2010).
\newblock The power of convex relaxation: Near-optimal matrix completion.
\newblock {\em IEEE Transactions on Information Theory}, 56(5):2053--2080.

\bibitem[Donoho, 2006]{donoho2006compressed}
Donoho, D.~L. (2006).
\newblock Compressed sensing.
\newblock {\em IEEE Transactions on information theory}, 52(4):1289--1306.

\bibitem[Fattahi and Sojoudi, 2020]{fattahi2020exact}
Fattahi, S. and Sojoudi, S. (2020).
\newblock Exact guarantees on the absence of spurious local minima for
  non-negative rank-1 robust principal component analysis.
\newblock {\em Journal of Machine Learning Research}, 21:1--51.

\bibitem[Ge et~al., 2017]{ge2017no}
Ge, R., Jin, C., and Zheng, Y. (2017).
\newblock No spurious local minima in nonconvex low rank problems: A unified
  geometric analysis.
\newblock In {\em Proceedings of the 34th International Conference on Machine
  Learning}, volume~70 of {\em Proceedings of Machine Learning Research}, pages
  1233--1242.

\bibitem[Jin et~al., 2021]{jin2021nonconvex}
Jin, C., Netrapalli, P., Ge, R., Kakade, S.~M., and Jordan, M.~I. (2021).
\newblock On nonconvex optimization for machine learning: Gradients,
  stochasticity, and saddle points.
\newblock {\em Journal of the ACM (JACM)}, 68(2):1--29.

\bibitem[Jin et~al., 2023]{jin2023understanding}
Jin, J., Li, Z., Lyu, K., Du, S.~S., and Lee, J.~D. (2023).
\newblock Understanding incremental learning of gradient descent: A
  fine-grained analysis of matrix sensing.
\newblock {\em arXiv preprint arXiv:2301.11500}.

\bibitem[Jin et~al., 2019]{jin2019towards}
Jin, M., Molybog, I., Mohammadi-Ghazi, R., and Lavaei, J. (2019).
\newblock Towards robust and scalable power system state estimation.
\newblock In {\em 2019 IEEE 58th Conference on Decision and Control (CDC)},
  pages 3245--3252. IEEE.

\bibitem[Li et~al., 2019]{li2019non}
Li, Q., Zhu, Z., and Tang, G. (2019).
\newblock The non-convex geometry of low-rank matrix optimization.
\newblock {\em Information and Inference: A Journal of the IMA}, 8(1):51--96.

\bibitem[Li et~al., 2018]{li2018algorithmic}
Li, Y., Ma, T., and Zhang, H. (2018).
\newblock Algorithmic regularization in over-parameterized matrix sensing and
  neural networks with quadratic activations.
\newblock In {\em Conference On Learning Theory}, pages 2--47. PMLR.

\bibitem[Ma and Fattahi, 2022]{ma2022global}
Ma, J. and Fattahi, S. (2022).
\newblock Global convergence of sub-gradient method for robust matrix recovery:
  Small initialization, noisy measurements, and over-parameterization.
\newblock {\em arXiv preprint arXiv:2202.08788}.

\bibitem[Ma et~al., 2022]{ma2022sharp}
Ma, Z., Bi, Y., Lavaei, J., and Sojoudi, S. (2022).
\newblock Sharp restricted isometry property bounds for low-rank matrix
  recovery problems with corrupted measurements.
\newblock In {\em Proceedings of the AAAI Conference on Artificial
  Intelligence}, volume~36, pages 7672--7681.

\bibitem[Ma et~al., 2023a]{ma2023geometric}
Ma, Z., Bi, Y., Lavaei, J., and Sojoudi, S. (2023a).
\newblock Geometric analysis of noisy low-rank matrix recovery in the exact
  parametrized and the overparametrized regimes.
\newblock {\em INFORMS Journal on Optimization}.

\bibitem[Ma et~al., 2023b]{ma2023over}
Ma, Z., Molybog, I., Lavaei, J., and Sojoudi, S. (2023b).
\newblock Over-parametrization via lifting for low-rank matrix sensing:
  Conversion of spurious solutions to strict saddle points.
\newblock In {\em International Conference on Machine Learning}. PMLR.

\bibitem[Ma and Sojoudi, 2023]{ma2023noisy}
Ma, Z. and Sojoudi, S. (2023).
\newblock Noisy low-rank matrix optimization: Geometry of local minima and
  convergence rate.
\newblock In {\em International Conference on Artificial Intelligence and
  Statistics}, pages 3125--3150. PMLR.

\bibitem[Molybog et~al., 2020]{molybog2020conic}
Molybog, I., Madani, R., and Lavaei, J. (2020).
\newblock Conic optimization for quadratic regression under sparse noise.
\newblock {\em The Journal of Machine Learning Research}, 21(1):7994--8029.

\bibitem[Recht et~al., 2010]{recht2010guaranteed}
Recht, B., Fazel, M., and Parrilo, P.~A. (2010).
\newblock Guaranteed minimum-rank solutions of linear matrix equations via
  nuclear norm minimization.
\newblock {\em SIAM Review}, 52(3):471--501.

\bibitem[Shechtman et~al., 2015]{shechtman2015phase}
Shechtman, Y., Eldar, Y.~C., Cohen, O., Chapman, H.~N., Miao, J., and Segev, M.
  (2015).
\newblock Phase retrieval with application to optical imaging: A contemporary
  overview.
\newblock {\em IEEE Signal Processing Magazine}, 32(3):87--109.

\bibitem[Singer, 2011]{singer2011angular}
Singer, A. (2011).
\newblock Angular synchronization by eigenvectors and semidefinite programming.
\newblock {\em Applied and Computational Harmonic Analysis}, 30(1):20--36.

\bibitem[St{\"o}ger and Soltanolkotabi, 2021]{stoger2021small}
St{\"o}ger, D. and Soltanolkotabi, M. (2021).
\newblock Small random initialization is akin to spectral learning:
  Optimization and generalization guarantees for overparameterized low-rank
  matrix reconstruction.
\newblock {\em Advances in Neural Information Processing Systems},
  34:23831--23843.

\bibitem[Wang et~al., 2017]{wang2017unified}
Wang, L., Zhang, X., and Gu, Q. (2017).
\newblock A unified computational and statistical framework for nonconvex
  low-rank matrix estimation.
\newblock In {\em Proceedings of the 20th International Conference on
  Artificial Intelligence and Statistics}, volume~54 of {\em Proceedings of
  Machine Learning Research}, pages 981--990.

\bibitem[Wright and Recht, 2022]{wright2022optimization}
Wright, S.~J. and Recht, B. (2022).
\newblock {\em Optimization for data analysis}.
\newblock Cambridge University Press.

\bibitem[Yalcin et~al., 2023]{yalcin2022semidefinite}
Yalcin, B., Ma, Z., Lavaei, J., and Sojoudi, S. (2023).
\newblock Semidefinite programming versus burer-monteiro factorization for
  matrix sensing.
\newblock In {\em Proceedings of the AAAI Conference on Artificial
  Intelligence}.

\bibitem[Zhang and Zhang, 2020]{zhang2020many}
Zhang, G. and Zhang, R.~Y. (2020).
\newblock How many samples is a good initial point worth in low-rank matrix
  recovery?
\newblock In {\em Advances in Neural Information Processing Systems},
  volume~33, pages 12583--12592.

\bibitem[Zhang et~al., 2021]{zhang2021general}
Zhang, H., Bi, Y., and Lavaei, J. (2021).
\newblock General low-rank matrix optimization: Geometric analysis and sharper
  bounds.
\newblock {\em Advances in Neural Information Processing Systems},
  34:27369--27380.

\bibitem[Zhang, 2022]{zhang2022improved}
Zhang, R.~Y. (2022).
\newblock Improved global guarantees for the nonconvex burer--monteiro
  factorization via rank overparameterization.
\newblock {\em arXiv preprint arXiv:2207.01789}.

\bibitem[Zhang et~al., 2019]{zhang2019sharp}
Zhang, R.~Y., Sojoudi, S., and Lavaei, J. (2019).
\newblock Sharp restricted isometry bounds for the inexistence of spurious
  local minima in nonconvex matrix recovery.
\newblock {\em Journal of Machine Learning Research}, 20(114):1--34.

\bibitem[Zhang et~al., 2017]{zhang2017conic}
Zhang, Y., Madani, R., and Lavaei, J. (2017).
\newblock Conic relaxations for power system state estimation with line
  measurements.
\newblock {\em IEEE Transactions on Control of Network Systems},
  5(3):1193--1205.

\bibitem[Zhu et~al., 2018]{zhu2018global}
Zhu, Z., Li, Q., Tang, G., and Wakin, M.~B. (2018).
\newblock Global optimality in low-rank matrix optimization.
\newblock {\em IEEE Transactions on Signal Processing}, 66(13):3614--3628.

\end{thebibliography}

\section*{Checklist}

The checklist follows the references. For each question, choose your answer from the three possible options: Yes, No, Not Applicable.  You are encouraged to include a justification to your answer, either by referencing the appropriate section of your paper or providing a brief inline description (1-2 sentences). 
Please do not modify the questions.  Note that the Checklist section does not count towards the page limit. Not including the checklist in the first submission won't result in desk rejection, although in such case we will ask you to upload it during the author response period and include it in camera ready (if accepted).

\textbf{In your paper, please delete this instructions block and only keep the Checklist section heading above along with the questions/answers below.}

 \begin{enumerate}

 \item For all models and algorithms presented, check if you include:
 \begin{enumerate}
   \item A clear description of the mathematical setting, assumptions, algorithm, and/or model. [Yes]
   \item An analysis of the properties and complexity (time, space, sample size) of any algorithm. [Yes]
   \item (Optional) Anonymized source code, with specification of all dependencies, including external libraries. [Yes]
 \end{enumerate}

 \item For any theoretical claim, check if you include:
 \begin{enumerate}
   \item Statements of the full set of assumptions of all theoretical results. [Yes]
   \item Complete proofs of all theoretical results. [Yes] Proofs are included in the Appendix (supplementary file).
   \item Clear explanations of any assumptions. [Yes]     
 \end{enumerate}

 \item For all figures and tables that present empirical results, check if you include:
 \begin{enumerate}
   \item The code, data, and instructions needed to reproduce the main experimental results (either in the supplemental material or as a URL). [Yes] An anonymous github link includes all the necessary code to reproduce the results.
   \item All the training details (e.g., data splits, hyperparameters, how they were chosen). [Yes]
         \item A clear definition of the specific measure or statistics and error bars (e.g., with respect to the random seed after running experiments multiple times). [Not Applicable]
         \item A description of the computing infrastructure used. (e.g., type of GPUs, internal cluster, or cloud provider). [Not Applicable] Since we're not interested in the running time of the experiments, any hardware should produce the same results.
 \end{enumerate}

 \item If you are using existing assets (e.g., code, data, models) or curating/releasing new assets, check if you include:
 \begin{enumerate}
   \item Citations of the creator If your work uses existing assets. [Not Applicable]
   \item The license information of the assets, if applicable. [Not Applicable]
   \item New assets either in the supplemental material or as a URL, if applicable. [Not Applicable]
   \item Information about consent from data providers/curators. [Not Applicable]
   \item Discussion of sensible content if applicable, e.g., personally identifiable information or offensive content. [Not Applicable]
 \end{enumerate}

 \item If you used crowdsourcing or conducted research with human subjects, check if you include:
 \begin{enumerate}
   \item The full text of instructions given to participants and screenshots. [Not Applicable]
   \item Descriptions of potential participant risks, with links to Institutional Review Board (IRB) approvals if applicable. [Not Applicable]
   \item The estimated hourly wage paid to participants and the total amount spent on participant compensation. [Not Applicable]
 \end{enumerate}

 \end{enumerate}

\appendix
\onecolumn
\newpage

\section{Appendix}

\subsection{Optimality Conditions}


\begin{lemma}\label{lem:lp_crit_cond}
	Given the problem \eqref{eq:l_p_problem}, its gradient and Hessian are given as
\begin{subequations}
\begin{align}
    &\langle \nabla f^l(X), U \rangle = \sum_{i=1}^m \langle A_i, XX^\top - M^* \rangle^{l-1} \langle A_i, UX^\top + XU^\top \rangle \quad \forall U \in \RR^{n \times r},
    \label{eq:focp_lp_f} 
\\
    &\nabla^2 f^l(X)[U,U] =  \sum_{i=1}^m 	\langle A_i, XX^\top - M^* \rangle^{l-2} \left[ (l-1) \langle A_i, UX^\top + XU^\top \rangle^2 + \langle A_i, XX^\top - M^* \rangle \langle A_i, 2 UU^\top \rangle \right] \quad \forall U \in \RR^{n \times r} 
    \label{eq:socp_lp_f}
\end{align}
\end{subequations}
\end{lemma}

\begin{lemma}\label{lem:lp_h_crit_cond}
	Given the problem \eqref{eq:lp_h_func}, its gradient, Hessian and high-order derivatives are equal to 
\begin{subequations}
 \begin{align}
&\nabla h^l(M) = \sum_{i=1}^m \langle A_i, M - M^* \rangle^{l-1} A_i, \label{eq:focp_lp} 
\\
&\nabla^2 h^l(M)[N,N] = (l-1) \sum_{i=1}^m 	\langle A_i, M - M^* \rangle^{l-2} \langle A_i, N \rangle^2 \quad \forall N \in \RR^{n \times n}, \label{eq:socp_lp}
  \\
&\nabla^p h^l(M)[\underbrace{N,\dots,N}_{\text{p times}}] = \frac{(l-1)!}{(l-p)!} \sum_{i=1}^m 	\langle A_i, M - M^* \rangle^{l-p} \langle A_i, N \rangle^p \quad \forall N \in \RR^{n \times n}
\end{align}
\end{subequations}
\end{lemma}

\subsection{Proofs in Section~\ref{sec:l2_case}}

\begin{proof}[Proof of Theorem~\ref{thm:l2_nospurious}]
Via the definition of RIP, we have that
\[
	h(M^*) \geq h(\hat X \hat X^\top) + \langle \nabla h(\hat X \hat X^\top), M^* - \hat X \hat X^\top \rangle + \frac{1-\delta}{2} \|\hat X \hat X^\top -M^* \|^2_F
 \]
 Given that $\hat X$ is a first-order critical point, it means that it must satisfy first-order optimality conditions, which means 
 \[
 	h(\hat X \hat X^\top) \hat X = 0,
 \]
 leading to 
 \[
 	\langle \nabla h(\hat X \hat X^\top), M^* \rangle \leq - \frac{1-\delta}{2} \|\hat X \hat X^\top -M^* \|^2_F
 \]
 since $h(\hat X \hat X^\top) - h(M^*) \geq 0$ via the construction of the objective function. Furthermore, as it is without loss of generality to assume the gradient of $h(M)$ to be symmetric \citep{zhang2021general}, and the fact that $M^*$ is positive-semidefinite, we have that 
 \[
		\langle \nabla h(\hat X \hat X^\top ), M^* \rangle \geq \lambda_{\min} (\nabla h(\hat X \hat X^\top )) \tr(M^*)
	\]
This leads to the fact that 
\begin{equation}\label{eq:minimum_eigenvalue_l2}
	\lambda_{\min} (\nabla h(\hat X \hat X^\top )) \leq -\frac{(1-\delta)  \|\hat X \hat X^\top -M^*\|^2_F}{2\tr(M^*)} \leq 0
\end{equation}
Given the optimality conditions, we know that for $\hat X$ to be a strict saddle, we must prove that there exists a direction $\Delta \in \RR^{n \times r}$ such that 
\begin{equation}\label{eq: hessian_saddle_cond_l2}
	2 \langle \nabla h (\hat X \hat X^\top) , \Delta \Delta^\top \rangle + \underbrace{\left[ \nabla^2 h(\hat X \hat X^\top) \right] (\hat X \Delta^\top + \Delta \hat X^\top, \hat X \Delta^\top + \Delta \hat X^\top)}_{P(\Delta)} < 0
\end{equation}
If we choose 
\[
	\Delta = uq^\top, \quad \|u\|_2, \|q\|_2 = 1, \quad \| \hat X q \|_2 = \sigma_r(\hat X), \  u^\top \nabla h(\hat X \hat X^\top ) u = \lambda_{\min} (\nabla h(\hat X \hat X^\top )),
\]
then it follows from the RIP condition that 
\begin{align*}
	P(\Delta) &\leq (1+\delta) \|\hat X \Delta^\top + \Delta \hat X^\top\|^2_F \\
	&= (1+\delta) \| u (\hat X q)^\top + (\hat X q) u^\top \|^2_F \\
	&= 2 (1+\delta) \|\hat X q\|^2_F + 2(1+\delta) \left( q^\top (\hat X^\top u) \right)^2 \\
	& = 2 (1+\delta) \sigma_r(\hat X)^2
\end{align*}
because of $\hat X^\top u = 0$ due to the first-order optimality condition. Therefore,
\begin{align*}
	2 \langle \nabla h (\hat X \hat X^\top) , \Delta \Delta^\top \rangle + P(\Delta) &\leq 2 \langle h (\hat X \hat X^\top) , \Delta \Delta^\top \rangle + 2 (1+\delta) \sigma_r(\hat X)^2 \\
	&= 2 \langle \nabla h (\hat X \hat X^\top) , uu^\top \rangle + 2 (1+\delta) \sigma_r(\hat X)^2 \\
	& = 2 u^\top \nabla h (\hat X \hat X^\top) u + 2 (1+\delta) \sigma_r(\hat X)^2 \\
	& = 2 \left( \lambda_{\min} (\nabla h(\hat X \hat X^\top )) + (1+\delta)  \sigma_r(\hat X)^2 \right) \\
	& \leq 2(1+\delta) \sigma_r(\hat X)^2 - \frac{(1-\delta)  \|\hat X \hat X^\top -M^*\|^2_F}{\tr(M^*)}
\end{align*}
Therefore, in order to make \eqref{eq: hessian_saddle_cond_l2} hold, we simply need
\[
	\|\hat X \hat X^\top -M^* \|^2_F > 2 \frac{1+\delta}{1-\delta} \tr(M^*) \sigma_r(\hat X)^2
\]
which concludes the proof.
\end{proof}

\begin{proof}[Proof for Theorem~\ref{thm:small_mstar}]
	An easy extension of Lemma B.2 from \cite{ma2023over} gives us that 
	\begin{equation}
		\sigma_r^2(\hat X) < \sqrt{\frac{2(1+\delta)}{r(1-\delta)}} \|M^*\|_F
	\end{equation}
	and we omit the proof for brevity. Then the only remaining step is to show that the right-hand side of \eqref{eq:close_criterion} is larger than that of \eqref{eq:l2_criterion}. This means the following equation
	\[
		(1-\delta^2) \lambda_{r^*}(M^*) \geq 2 \frac{1+\delta}{1-\delta} \tr(M^*) \sqrt{\frac{2(1+\delta)}{r(1-\delta)}} \|M^*\|_F
	\]
	is sufficient to 
	\[
		(1-\delta^2) \lambda_{r^*}(M^*) \geq 2 \frac{1+\delta}{1-\delta} \tr(M^*) \sigma_r^2(\hat X)
	\]
	which yields \eqref{eq:mstar_criterion} after rearrangements.
\end{proof}

\subsection{Proofs in Section~\ref{sec:modified_loss}}

Before proceeding to the main proof, we first present a technical lemma:
\begin{lemma}\label{lem:lp_norm_inequality}
	Given a vector $x \in \RR^n$,  we have that 
	\begin{equation}
		\|x\|_p \leq n^{1/p-1/q} \|x\|_q \quad \forall  \ q \geq p
	\end{equation}
\end{lemma}

\begin{proof}[Proof of Lemma~\ref{lem:lp_norm_inequality}]
	By applying Holder's inequality, we obtain that 
\[
	\sum\limits_{i=1}^n |x_i|^p=
\sum\limits_{i=1}^n |x_i|^p\cdot 1\leq
\left(\sum\limits_{i=1}^n (|x_i|^p)^{\frac{q}{p}}\right)^{\frac{p}{q}}
\left(\sum\limits_{i=1}^n 1^{\frac{q}{q-p}}\right)^{1-\frac{p}{q}}=
\left(\sum\limits_{i=1}^n |x_i|^q\right)^{\frac{p}{q}} n^{1-\frac{p}{q}}
\]
Then
\[
	\Vert x\Vert_p=
\left(\sum\limits_{i=1}^n |x_i|^p\right)^{1/p}\leq
\left(\left(\sum\limits_{i=1}^n |x_i|^q\right)^{\frac{p}{q}} n^{1-\frac{p}{q}}\right)^{1/p}=
\left(\sum\limits_{i=1}^n |x_i|^q\right)^{\frac{1}{q}} n^{\frac{1}{p}-\frac{1}{q}}=
n^{1/p-1/q}\Vert x\Vert_q
\]
	
\end{proof}

\begin{proof}[Proof to Theorem~\ref{thm:lp_nospurious}]
First, we define
\[
	h_\lambda^l(M) = h(M) + \lambda h^l(M)
\]
	Then, focusing on $h^l(\cdot)$, using Taylor's theorem and the mean-value theorem, we have
	\begin{equation}
		\begin{aligned}
			h^l(M^*) = &h^l(\hat X \hat X^\top) + \langle \nabla h^l (\hat X \hat X^\top), \Delta \rangle + \frac{1}{2!} \nabla^2 h^l(\hat X \hat X^\top) \left[\Delta,\Delta \right] + \\
			& \frac{1}{3!}\nabla^3 h^l(\hat X \hat X^\top) \left[\Delta,\Delta, \Delta \right] + \dots + \frac{1}{l!}\nabla^l h^l(\tilde M) \left[\Delta,\dots,\Delta \right]
		\end{aligned}
	\end{equation}
	where $\tilde M$ is a convex combination of $M^*$ and $\hat X \hat X^\top$ and
	\[
		\Delta = M^* - \hat X \hat X^\top .
	\]
	Using Lemma~\ref{lem:lp_h_crit_cond}, we know that 
	\begin{align*}
		&\frac{1}{p!} \nabla^p h^l(\hat X \hat X^\top) \left[\Delta,\dots,\Delta \right] = \frac{1}{p!} \sum_{i=1}^m 	\langle A_i, \Delta \rangle^{l} \quad \forall p \in [2,l-1], \\
		&\frac{1}{l!}\nabla^l h^l(M) \left[\Delta,\dots,\Delta \right] = \frac{1}{l!} \sum_{i=1}^m 	\langle A_i, \Delta \rangle^{l} \quad \forall M \in \RR^{n \times n}
	\end{align*}
	Hence, it is possible to rewrite
	\begin{equation}\label{eq:lp_taylors}
		\begin{aligned}
			h^l(M^*) &=h^l(\hat X \hat X^\top) + \langle \nabla h^l (\hat X \hat X^\top), \Delta \rangle + \sum_{p=2}^l \frac{(l-1)!}{(l-p)!p!} \sum_{i=1}^m 	\langle A_i, \Delta \rangle^{l} \\
			&=h^l(\hat X \hat X^\top) + \langle \nabla h^l (\hat X \hat X^\top), \Delta \rangle + \sum_{i=1}^m \langle A_i, \Delta \rangle^{l} \sum_{p=2}^l \frac{(l-1)!}{(l-p)!p!}  	 \\
			&=h^l(\hat X \hat X^\top) + \langle \nabla h^l (\hat X \hat X^\top), M^*-\hat X \hat X^\top  \rangle + (\frac{2^l-1}{l}-1)\|\mathcal{A}(\hat X \hat X^\top -M^*)\|^l_l
		\end{aligned}
	\end{equation}
	By Lemma~\ref{lem:lp_norm_inequality}, if $l > 2$, it holds that 
	\[
		\|\mathcal{A}(\hat X \hat X^\top -M^*)\|^l_l \geq \|\mathcal{A}(\hat X \hat X^\top -M^*)\|^l_2 m^{(2-l)/2}
	\]
	Furthermore, combining this with the RIP property gives rise to
	\begin{equation}
		\|\mathcal{A}(\hat X \hat X^\top -M^*)\|^l_l \geq (1-\delta)^{l/2} \|\hat X \hat X^\top - M^*\|_F^l m^{(2-l)/2}
	\end{equation}
	Therefore we know that 
	\begin{equation}
	\begin{aligned}
		h_\lambda^l(M^*) = h(M^*) + \lambda h^l(M^*) \geq h_\lambda^l(\hat X \hat X^\top) + \langle \nabla h_\lambda^l(\hat X \hat X^\top), M^* - \hat X \hat X^\top \rangle + L
	\end{aligned}
	\end{equation}
	in which
	\[
		L \coloneqq \frac{1-\delta}{2} \|M^*-\hat X \hat X^\top \|^2_F + \lambda (1-\delta)^{l/2} C(l) \|M^*-\hat X \hat X^\top \|^l_F
	\]
	If $\hat X$ is a first-order critical point, a repeated application of \eqref{eq:focp_lp_f} yields that
	\[
		\nabla h_\lambda^l (\hat X \hat X^\top) \hat X = 0 \implies 
		\langle \nabla h_\lambda^l (\hat X \hat X^\top), \hat X \hat X^\top \rangle = 0
	\]
	which after rearrangement leads to 
	\begin{equation}
	\begin{aligned}
		\langle \nabla h_\lambda^l (\hat X \hat X^\top), M^* \rangle &\leq \left[ h_\lambda^l(M^*)-h_\lambda^l(\hat X \hat X^\top)\right] - L \\
		&\leq - L
	\end{aligned}
	\end{equation}
	where the second inequality follows from the fact that $M^*$ is the global minimizer of \eqref{eq:l_p_problem}. Since the sensing matrices can be assumed to be symmetric without loss of generality \citep{zhang2021general}, $\nabla h^l (\hat X \hat X^\top)$ can be assumed to be symmetric according to \eqref{eq:focp_lp}. This means that 
	\[
		\langle \nabla h_\lambda^l (\hat X \hat X^\top), M^* \rangle \geq \tr(M^*) \lambda_{\min}(\nabla h_\lambda^l (\hat X \hat X^\top))
	\]
	which further leads to
	\begin{equation}\label{eq:fp_hessian_lambdmin}
		\lambda_{\min} (\nabla h_\lambda^l (\hat X \hat X^\top)) \leq - \frac{L }{\tr(M^*)}
	\end{equation}
	Now, we turn to the Hessian of $f_\lambda^l(\cdot)$, which given \eqref{eq:socp_lp_f} is
	\begin{equation}\label{eq:f_socp_BC}
	\begin{aligned}
		\nabla^2 f_\lambda^l(\hat X)[U,U] =  \underbrace{\sum_{i=1}^m \left[\lambda (l-1)	\langle A_i, \hat X \hat X^\top - M^* \rangle^{l-2} + 1 \right]   \langle A_i, U \hat X^\top + \hat XU^\top \rangle^2 }_{B} \\
		+  \underbrace{ \sum_{i=1}^m \left[ \lambda \langle A_i, \hat X \hat X^\top - M^* \rangle^{l-1} + \langle A_i, \hat X \hat X^\top - M^* \rangle \right] \langle A_i, 2 UU^\top \rangle}_{A}  \quad \forall U \in \RR^{n \times r} 
	\end{aligned}
	\end{equation}
	Since
	\[
		\langle A_i, U \hat X^\top + \hat XU^\top \rangle^2 \leq \sum_{i=1}^m \langle A_i, U \hat X^\top + \hat XU^\top \rangle^2 \leq (1+\delta) \| U\hat X^\top + \hat XU^\top\|^2_F \quad \forall i
	\]
	then if we choose $U$ such that
	\[
	U = uq^\top, \quad \|u\|_2, \|q\|_2 = 1, \quad \| \hat X q \|_2 = \sigma_r(\hat X), \  u^\top \nabla h^l(\hat X \hat X^\top ) u = \lambda_{\min} (\nabla h_\lambda^l(\hat X \hat X^\top )),
\]
it can be shown that
\begin{align*}
	(1+\delta) \| U\hat X^\top + \hat XU^\top\|^2_F &= (1+\delta) \| u (\hat X q)^\top + (\hat X q) u^\top \|^2_F \\
	&= 2 (1+\delta) \|\hat X q\|^2_F + 2(1+\delta) \left( q^\top (\hat X^\top u) \right)^2 \\
	& = 2 (1+\delta) \sigma_r(\hat X)^2
\end{align*}
since $(\hat X q) u^\top = 0$ as $u$ is an eigenvector of $\nabla h_\lambda^l(\hat X \hat X^\top )$, which is orthogonal to $\hat X$ as required by \eqref{eq:focp_lp_f}. This further implies that 
\begin{equation} \label{eq:B_derive}
	\begin{aligned}
		B &\leq 2(1+\delta) \sigma_r(\hat X)^2
 \left[ \lambda (l-1) \sum_{i=1}^m \langle A_i, \hat X \hat X^\top - M^* \rangle^{l-2} + 1 \right] \\
 &= 2(1+\delta) \sigma_r(\hat X)^2
 \left[ \lambda (l-1) \| \mathcal{A}(\hat X \hat X^\top - M^*) \|^{l-2}_{l-2} + 1 \right]  \\
 &\leq 2(1+\delta) \sigma_r(\hat X)^2
 \left[ \lambda (l-1) (\| \mathcal{A}(\hat X \hat X^\top - M^*) \|_2 )^{l-2} + 1 \right] \\
 &\leq 2(1+\delta) \sigma_r(\hat X)^2
 \left[ \lambda (l-1) (\sqrt{1+\delta} \|\hat X \hat X^\top - M^*\|_F )^{l-2} + 1 \right]  \\
 & = 2(1+\delta) \sigma_r(\hat X)^2
 \left[ \lambda (l-1) (1+\delta)^{(l-2)/2} \|\hat X \hat X^\top - M^*\|_F^{l-2} + 1 \right]
	\end{aligned}
\end{equation}
Also, given \eqref{eq:focp_lp} and our choice of $U$, it is apparent that 
\begin{equation}\label{eq:C_derive}
	A = 2 u^\top \nabla h_\lambda^l(\hat X \hat X^\top) u = 2 \lambda_{\min} (\nabla h_\lambda^l(\hat X \hat X^\top))
\end{equation}
Therefore, given \eqref{eq:B_derive}, \eqref{eq:C_derive} and \eqref{eq:fp_hessian_lambdmin}, by substituting them back into \eqref{eq:f_socp_BC}, we arrive at 
\begin{equation}
	\nabla^2 f^l(\hat X)[U,U] \leq 2 \sigma_r(\hat X)^2 \left( \lambda (l-1) (1+\delta)^{l/2}\|\hat X \hat X^\top - M^*\|_F^{l-2} + (1+\delta) \right) - 2L/\tr(M^*)
\end{equation}
so the right-hand side of the above equation is strictly negative if \eqref{eq:lp_criterion} is satisfied.
\end{proof}

\end{document}